\documentclass[twocolumn]{autart}    

\usepackage{graphicx}          
\usepackage{amsmath}
\usepackage{amsfonts}
\usepackage{amssymb}
\usepackage{mathrsfs}
\usepackage{multirow}
\usepackage{enumerate}
\usepackage{float}
\usepackage[caption = false]{subfig}
\usepackage{mathtools, cuted}
\usepackage{xcolor}

\def \blkdiag {\operatorname{block-diag}}

\newtheorem{theorem}{Theorem}
\newtheorem{lemma}{Lemma}

\newtheorem{remark}{Remark}

\theoremstyle{definition}
\newtheorem{definition}[theorem]{Definition}

\begin{document}

\begin{frontmatter}

\title{A Double-Layer Jacobi Method for PDE-Constrained Nonlinear Model Predictive Control\thanksref{footnoteinfo}} 

\thanks[footnoteinfo]{This work was partly supported by JSPS KAKENHI Grant Number 15H02257. Corresponding author H. Deng. Tel. +81-80-6470-1087. Fax. +81-75-753-5042.}

\author[Kyoto]{Haoyang Deng}\ead{deng.haoyang.23r@st.kyoto-u.ac.jp},    
\author[Kyoto]{Toshiyuki Ohtsuka}\ead{ohtsuka@i.kyoto-u.ac.jp}               

\address[Kyoto]{Department of Systems Science, Graduate School of Informatics, Kyoto University, Sakyo-ku, Kyoto  606-8501, Japan}  


\begin{abstract}                          
	This paper presents a real-time optimization method for nonlinear model predictive control (NMPC) of systems governed by partial differential equations (PDEs). 
	The NMPC problem to be solved is formulated by discretizing the PDE system in space and time by using the finite difference method. 
	The proposed method is called the double-layer Jacobi method, which exploits both the spatial and temporal sparsities of the PDE-constrained NMPC problem. 
	In the upper layer, the NMPC problem is solved by ignoring the temporal couplings of either the state or costate (Lagrange multiplier corresponding to the state equation) equations so that the spatial sparsity is preserved. 
	The lower-layer Jacobi method is a linear solver dedicated to PDE-constrained NMPC problems by exploiting the spatial sparsity. 
	Convergence analysis indicates that the convergence of the proposed method is related to the prediction horizon and regularization.
	Results of a numerical experiment of controlling a heat transfer process show that the proposed method is two orders of magnitude faster than the conventional structure-exploiting Newton's method.
\end{abstract}

\end{frontmatter}

\section{Introduction}
Nonlinear model predictive control (NMPC), also referred to as nonlinear receding horizon control or nonlinear moving horizon control, is an optimization-based control method for nonlinear systems. 
The control input is obtained by solving an optimization problem that usually minimizes a tracking cost under the constraints of the system dynamics. 
Moreover, general input and state constraints or economic costs can be integrated into the optimization problem, making NMPC a powerful advanced control technique and a popular research topic. 
On the other hand, the generality of NMPC brings about computational difficulties in solving the underlying optimization problem in real time. 
Considerable efforts and progress have been made toward the real-time NMPC control of systems described by ordinary differential equations (ODEs) in recent years, ranging over automatic code generation \cite{ohtsuka2002automatic}, first-order iteration (e.g., \cite{englert2019software}), and parallel computing \cite{deng2019parallel}, to name just a few. 

Besides the NMPC control of ODE systems, NMPC control of systems described by partial differential equations (PDEs), such as Navier-Stokes equations for fluid flow and heat transfer equations for chemical processes, has gained increasing attention due to the optimal and constraint-handling properties of NMPC. 
However, PDE-constrained NMPC presents a great challenge for real-time optimization due to the infinite-dimensional state space of PDE systems. 
The solution methods for PDE-constrained NMPC can be generally categorized as indirect or direct. 
Indirect methods analytically derive the optimality conditions for PDE-constrained NMPC and then solve these conditions numerically. 
For example, in \cite{hashimoto2012receding}, the analytic optimality conditions for the NMPC control of a class of parabolic PDEs are derived. 
These optimality conditions are discretized to a set of  nonlinear algebraic equations, the solution of which is then traced by using the C/GMRES  method \cite{ohtsuka2004continuation}. 
Moreover, the so-called contraction mapping method can be applied efficiently if the nonlinear algebraic equations satisfy certain structure conditions. 
In contrast to indirect methods, which ``first optimize, then discretize," direct methods need PDE systems to be first discretized both in space and time. 
The NMPC problem is formulated on the basis of the discretized PDE system, which leads to a nonlinear program (NLP). 
Since a fine-grained spatial discretization results in a large number of states or optimization variables, model reduction techniques are frequently applied. 
A common way, e.g., \cite{ou2009model}, is to combine the proper orthogonal decomposition method \cite{sirovich1987turbulence} with Galerkin projection to obtain a low-dimensional dynamical system. 
An alternative approach is to use Koopman operator-based reduced order models \cite{peitz2019koopman}. 
The key idea is to transform the dynamical system into switched autonomous systems by restricting the control input to a finite number of constant values.
The switched autonomous systems are approximated by low-dimensional linear systems using the Koopman operator. 

Unlike the model reduction methods, this paper deals directly with the discretized PDE system. 
Although conventional NMPC methods for ODE systems can in principle be applied, they are computationally expensive due to the large number of states. 
For example, structure-exploiting methods for NMPC, e.g., \cite{steinbach1994structured} and \cite{zanelli2020forces}, perform Riccati recursions and have computational complexities of $ \mathcal{O}(N(n_u+n_x)^3) $, where $ N $ is the number of the temporal discretization grid points and $ n_u $ and $ n_x $ are the numbers of the control inputs and states, respectively. 
Even the state variables can be first eliminated in condensing-based methods, the elimination procedure is roughly of $ \mathcal{O}(N^2n_x^2) $ \cite{andersson2013condensing}. 
Note that sparsity exists both in the spatial and temporal directions of PDE-constrained NMPC problems. 
Structure-exploiting methods can only make use of the temporal sparsity along the prediction horizon, and the spatial sparsity is destroyed due to Riccati recursion. 

In this paper, we present a double-layer Jacobi method that exploits both the temporal and spatial sparsities. 
The upper-layer Jacobi (or Jacobi-type) method is derived by ignoring the temporal couplings of either the state or costate equations such that the spatial sparsity can be preserved. 
The upper-layer Jacobi method can be applied to the NMPC control of not only PDE systems, but also general nonlinear systems. 
Convergence analysis shows that the upper-layer Jacobi method can be guaranteed to converge by introducing regularization and choosing a short prediction horizon. 
The lower-layer Jacobi method is an iterative linear solver tailored to PDE-constrained NMPC by exploiting the spatial sparsity.
The performance of the proposed method is assessed by controlling the temperature distribution of a two-dimensional heat transfer process on a thin plate. 
The proposed method is matrix-free and has a complexity of $ \mathcal{O}(N(n_u+n_x)) $ for the heat transfer example. The numerical example shows that the proposed method is two orders of magnitude faster than the conventional structure-exploiting method. 

This paper is organized as follows. 
PDE systems and their spatial discretization are described in Section \ref{sec_pde}.
The NMPC problem and its Karush-Kuhn-Tucker (KKT) conditions are given in Section \ref{sec_nmpc}. 
The proposed double-layer Jacobi method is introduced in Section \ref{sec_doublelayer}. 
Section \ref{sec_exp} demonstrates the performance of the proposed method. 
Finally, this paper is summarized in Section \ref{sec_conclusion}.

\subsection{Notations}
Let $ v_{(i)} $ be the $ i $-th component of a vector $ v\in\mathbb{R}^{n} $. 
For a matrix $ P\in\mathbb{R}^{n\times n} $, we denote $ \rho(P) $ as the spectral radius of $ P $. 
The symbol $ \|\cdot\| $ denotes the Euclidean norm for a vector and the Frobenius norm for a matrix. 
The weighted norm is defined as $ \|v\|_P := \sqrt{v^TPv}$. 
For an iteration variable $ s $, we denote $ s^k $ as the value of $ s $ at the $ k $-th iteration and $ s^* $ as the optimal solution or fixed point. 
For a differentiable function $ f(v): \mathbb{R}^{n} \rightarrow \mathbb{R}^{m}$, we denote $ \nabla_{v}f \in \mathbb{R}^{m\times n}$ as the Jacobian matrix of $ f $. 
Identity and zero matrices are denoted by $ I $ and $ 0 $, respectively, and their sizes are indicated by using subscripts if necessary. 

\section{PDE systems}\label{sec_pde}

In this paper, we consider the NMPC control of a general class of PDE systems defined on the spatial domain $ \Omega\subset  \mathbb{R}^{n}$ and temporal domain $ \Gamma \subset\mathbb{R}$:
\begin{equation}\label{eq_pde}
\begin{split}
a&(u(t),w(p,t)) \frac{\partial^2 w(p,t)}{\partial t^2}+ b(u(t),w(p,t))\frac{\partial w(p,t)}{\partial t}\\
&= c(u(t),w(p,t)) \triangle w(p,t) + d(u(t),w(p,t)),
\end{split}
\end{equation}
where $ u\in \mathbb{U}\subset \mathbb{R}^{n_u} $ is the control input, $ w \in \mathbb{W}\subset \mathbb{R}$ is the PDE state, $ \triangle w(p,t):= \sum_{i=1}^{n}{\partial^2 w(p,t)}/{\partial p_{(i)}^2}$ denotes the Laplacian of $ w $, and $ a $, $ b $, $ c$ and $ d: \mathbb{U}\times \mathbb{W}\rightarrow \mathbb{R}$ are twice-differentiable nonlinear functions of $ u $ and $ w $. 
The boundary conditions, such as the Dirichlet and Neumann boundary conditions,  can be given to be input- and state-dependent, i.e., as functions of $ u $ and $ w $.
We assume that $ a $ is not zero for every $ (u,w) \in \mathbb{U} \times \mathbb{W} $ when the second-order time derivative of $ w $ is involved.

Note that \eqref{eq_pde} is a very general description of PDE systems. 
Many of the PDE systems, such as the heat transfer equation and wave equation, fall into the form of \eqref{eq_pde}. 
For PDE systems that are not in this form, e.g., the Navier-Stokes equations including both gradients and algebraic variables, we will discuss later in Remark \ref{rmk_second_general} that the results of this paper can in principle be extended. 

\subsection{Spatial discretization}
\begin{figure}
	\begin{center}
		\includegraphics[width=0.4\textwidth]{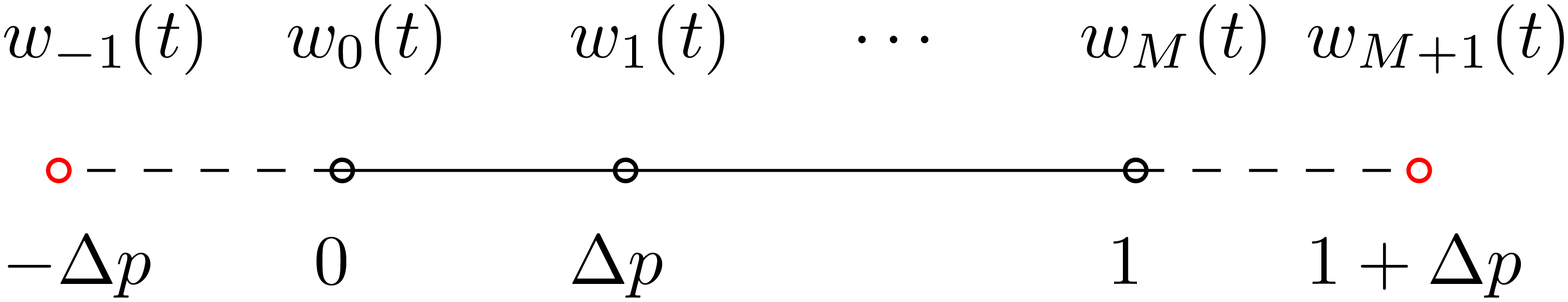}
	\end{center}
	\caption{Spatial discretization points and fictitious points}\label{fig_1ddis}
\end{figure}

We first introduce the spatial discretization of \eqref{eq_pde} by using the finite difference method. 
Without loss of generality, we demonstrate the discretization by using a one-dimensional system on an unit space interval $ \Omega:=[0,1] $ satisfying the following Neumann boundary condition:
\begin{equation}\label{eq_neumann}
\frac{\partial w(p,t)}{\partial p} = e(u(t),w(p,t)),\ p=0\ \text{and}\ 1.
\end{equation}
Let $ M+1 $ be the number of the spatial discretization grid points and $ \Delta p $ be the corresponding step size. 
The finite difference method is to approximate derivatives by using finite differences, i.e., for $  j\in\{0,\cdots,M\} $,
\begin{equation*}
\begin{split}
\triangle w(j\Delta p,t) \approx 
\frac{w_{j+1}(t) - 2w_{j}(t) + w_{j-1}(t) }{\Delta p^2},
\end{split}
\end{equation*}
where $ w_j(t) := w(j\Delta p,t) $. At $ j=0 $ and $ M $, two fictitious points $ w_{-1}(t) $ and $ w_{M+1}(t)$, as illustrated in Fig. \ref{fig_1ddis}, are introduced to deal with the boundary condition. 
By using the finite difference method, the Neumann boundary condition \eqref{eq_neumann} translates into the difference equations as follows.
\begin{subequations}\label{eq_dis_neumann}
	\begin{align}
	\frac{w_1(t)-w_{-1}(t)}{2\Delta p} &= e(u(t),w_0(t))\\
		\frac{w_{M+1}(t)-w_{M-1}(t)}{2\Delta p} &= e(u(t),w_M(t))
	\end{align}
\end{subequations}
The PDE system \eqref{eq_pde} is then discretized into
\begin{equation}\label{eq_pde_des}
\begin{split}
a&(u(t),w_j(t))\ddot{w}_j(t) + b(u(t),w_j(t))\dot{w}_j(t) \\
&= c(u(t),w_j(t))\frac{w_{j+1}(t)  - 2w_{j}(t) + w_{j-1}(t)}{\Delta p^2} \\
&\quad + d(u(t),w_j(t)), \ j\in\{0,\cdots,M\},
\end{split} 
\end{equation}
where the fictitious points $ w_{-1}(t) $ and $ w_{M+1}(t)$ can be eliminated by using the discretized boundary conditions \eqref{eq_dis_neumann}.
Note that the discretized PDE system \eqref{eq_pde_des} is described by a finite number of states:
\begin{equation*}
\begin{split}
x(t) &:=(W(t),\dot{W}(t))\\
     &:= (w_0(t),\cdots,w_M(t),\dot{w}_0(t),\cdots,\dot{w}_M(t))\in\mathbb{R}^{n_x},
\end{split}
\end{equation*}
where $ n_x=2(M+1) $. 
The dynamics of the discretized PDE system are given by 
\begin{equation}\label{eq_ode}
\dot{x}(t) = 
\left[
\begin{array}{c}
\dot{W}(t)\\ 
g(u(t),x(t))
\end{array} 
\right]=:
f(u(t),x(t)),
\end{equation}
where $ g(u(t),x(t)) $ denotes the expression of $ \ddot{W}(t) $ obtained from \eqref{eq_pde_des}. 

\section{NMPC}\label{sec_nmpc}
The spatial discretization approximates the PDE system \eqref{eq_pde} into an ODE system \eqref{eq_ode}, which is used to design the NMPC controller. 
For a prediction horizon $ T>0 $, we consider the following $ N $-stage NMPC problem based on the backward Euler's method with a temporal discretization step size of $ h:=T/N $:
\begin{equation}\label{eq:nlp}
\begin{split}
&\min_{\substack{u_1,\cdots,u_N, \\ x_0,\cdots,x_N}}\sum_{i=1}^{N}hl_i(u_i,x_i)\\
&\begin{split}
\text{s.t.} \quad 		& x_0 = \bar{x}_0,\\
& x_i = x_{i-1} + hf(u_i,x_i), \quad i\in\{1,\cdots,N\},\\
& u_i \in \mathbb{U},\ x_i\in \mathbb{X}, \quad i\in\{1,\cdots,N\},
\end{split}
\end{split}
\end{equation}
where $ \bar{x}_0 $ is the initial state, $ \mathbb{U} := [\underline{u},\bar{u}] $ and $ \mathbb{X} := [\underline{x},\bar{x}]$ are the admissible sets of $ u $ and $ x $ (the boundaries $ \underline{u} $, $ \bar{u} $, $ \underline{x} $, and $ \bar{x} $ are given), respectively, and $ l_i(u,x): \mathbb{U}\times \mathbb{X} \rightarrow \mathbb{R}$ is the stage cost function, which is assumed to be twice differentiable. 

\subsection{Regularization and relaxation}
Instead of solving the original NMPC problem \eqref{eq:nlp}, we add the following regularization term to the cost:
\begin{equation}\label{eq_reg}
\frac{\gamma}{2}\|u_i-\tilde{u}_i^*\|^2,
\end{equation}
where $ \gamma\geq 0 $ is the regularization parameter and $ \tilde{u}_i^* $ is the given regularization reference  and regarded as the estimation of the optimal control input of \eqref{eq:nlp}. 
We discuss later in Remark \ref{rmk_gamma} the role of regularization and the selection of $ \tilde{u}_i^* $. 

We adopt the interior-point method to relax the regularized NMPC problem by transferring the inequality constraints into a logarithmic barrier function added to the cost. 
To simplify the notations, the inequality constraints are put into a single vector-valued function 
\begin{equation*}
G(u,x)\geq 0.
\end{equation*}
We obtain the following relaxed regularized NMPC (RR-NMPC) problem:
\begin{equation}\label{eq:nlp_relaxed}
\begin{split}
&\min_{\substack{u_1,\cdots,u_N, \\ x_0,\cdots,x_N}}\sum_{i=1}^{N}\left(hl_i(u_i,x_i)+h \Phi(u,x)+\frac{\gamma}{2}\|u_i-\tilde{u}_i^*\|^2\right)\\
&\begin{split}
\text{s.t.} \quad 		& x_0 = \bar{x}_0,\\
& x_i = x_{i-1} + hf(u_i,x_i), \quad i\in\{1,\cdots,N\},
\end{split}
\end{split}
\end{equation}
where 
$ \Phi(u,x) := -\tau\sum_{j}\ln G_{(j)}(u,x)$ ($ \tau>0 $ is the barrier parameter). The RR-NMPC problem \eqref{eq:nlp_relaxed} approaches the original NMPC problem \eqref{eq:nlp} when $ \tau\rightarrow 0 $ and either $ \gamma =0 $ or $ \tilde{u}_i^* $ is given to be the optimal control input of \eqref{eq:nlp}.

For the consideration of feasibility, the state constraints usually need to be softened. The softened NMPC problem can also be regularized, relaxed, and written in the form of \eqref{eq:nlp_relaxed}. 
\subsection{KKT conditions}
Let $ \lambda_i\in\mathbb{R}^{n_x} $ be the Lagrange multiplier  (costate)  corresponding to the $ i $-th state equation. 
For the sake of brevity, we define
\begin{equation*}
s :=  (x,u,\lambda)\ \text{and}\ S:= (s_1,\cdots,s_N).
\end{equation*}
Let $ \mathcal{H}_i(s)$ be the Hamiltonian defined by
\begin{equation*}
\mathcal{H}_i(s) := l_i(u,x) + \Phi(u,x) + \lambda^Tf(u,x).
\end{equation*}
Let $ \mathcal{K}_i(x_{i-1},s_i,\lambda_{i+1}) $ be defined by
\begin{equation*}
\begin{split}
\mathcal{K}_i(x_{i-1},s_i,\lambda_{i+1})
:=
\left[
\begin{array}{l}
x_{i-1}- x_i+ hf(u_i,x_i) \\ 
h\nabla_u\mathcal{H}_i(s_i)^{T} + \gamma (u_i-\tilde{u}_i^*)\\
\lambda_{i+1}- \lambda_i+ h \nabla_x\mathcal{H}_i(s_i)^{T}
\end{array}
\right]
\end{split}
\end{equation*}
with $ x_0 = \bar{x}_0 $ and $ \lambda_{N+1}=0 $.
The KKT conditions for the RR-NMPC problem \eqref{eq:nlp_relaxed} are 
\begin{equation}\label{eq_KKT}
\mathcal{K}_i(x_{i-1}^*,s_i^*,\lambda_{i+1}^*)=0,\  i\in\{1,\cdots,N\}. 
\end{equation}
Although the KKT conditions \eqref{eq_KKT} are only the necessary conditions for optimality, we solve the RR-NMPC problem \eqref{eq:nlp_relaxed} by solving the nonlinear algebraic equations \eqref{eq_KKT}. 

We introduce the following shorthand at the $ k $-th iteration:
\begin{equation*}
\begin{split}
&\mathcal{K}_i^{k} := \mathcal{K}_i(x_{i-1}^k,s_i^k,\lambda_{i+1}^k),\  \mathcal{K}^{k}:=(\mathcal{K}_1^{k},\cdots,\mathcal{K}_N^{k}), \\ &\nabla_{s_i}\mathcal{K}_i^k:=\nabla_{s_i}\mathcal{K}_i(x_{i-1}^k,s_i^k,\lambda_{i+1}^k), \quad \text{etc}.
\end{split}
\end{equation*}

\subsection{Newton's method}\label{sec_reg_newton}
In solving the KKT conditions \eqref{eq_KKT} by using the Newton's method, the search direction $ \Delta S^k := (\Delta s_1^k,\cdots, \Delta s_N^k) $ is obtained by solving the following KKT system:
\begin{equation}\label{KKT_sys}
\left[
\begin{array}{lllll}
\ddots &      &      &      &       \\
\cdots & D^k_{i-1}    & M_U &      &       \\
& M_L & D^k_{i}     & M_U &       \\
&      & M_L & D^k_{i+1}     & \cdots \\
&      &      &      & \ddots
\end{array}
\right]
\left[
\begin{array}{l}
\ \ \ \vdots\\
\Delta s_{i-1}^k\\ 
\Delta s_i^k\\ 
\Delta s_{i+1}^k\\ 
\ \ \ \vdots
\end{array} 
\right]=
\left[
\begin{array}{l}
\ \ \ \vdots\\ 
\mathcal{K}_{i-1}^{k}\\
\mathcal{K}_i^{k}\\
\mathcal{K}_{i+1}^{k}\\
\ \ \ \vdots
\end{array} 
\right].
\end{equation}
Here, 
$ D_i^k:= \nabla_{s_i}\mathcal{K}_i^{k} $ (the expression can be found in \eqref{eq_decom}) and the constant matrices $ M_L $ and $ M_U $ given as follows show the couplings of the state and costate equations, respectively.  
\begin{equation*}
M_L := 
\left[
\begin{array}{ccc}
I_{n_x} & 0 & 0 \\
0& 0  &0  \\
0&0 & 0
\end{array}
\right] 
\ 
M_U := 
\left[
\begin{array}{ccc}
0 & 0 & 0 \\
0& 0  &0  \\
0&0 & I_{n_x}
\end{array}
\right]
\end{equation*}
After the search direction is calculated, a line search is performed to guarantee the primal feasibility $ G(u,x) \geq 0$, i.e.,
\begin{equation*}
S^{k+1} = S^{k} - \alpha^{\max}\Delta S^k, \  i\in\{1,\cdots,N\},
\end{equation*}
where $ \alpha^{\max} $ is obtained from the fraction-to-the-boundary rule \cite{wright1999numerical}:
\begin{equation}\label{eq_alphamax}
\begin{split}
\alpha^{\max} = \max\{\alpha^{\max}\in (0,1]:\ &G_i^{k+1}\geq 0.005 G_i^k,\\&\ i\in\{1,\cdots,N\} \}.
\end{split}
\end{equation}
Since the KKT matrix in \eqref{KKT_sys} is a block-tridiagonal matrix, solving \eqref{KKT_sys} by using the block Gaussian elimination method has a computational complexity of $ \mathcal{O}(N(n_x+n_u)^3) $. 
Many of the existing structure-exploiting methods tailored to NMPC are of this complexity. 
Note that sparsity exists in both the upper-level KKT matrix in \eqref{KKT_sys}  and the lower-level Jacobian matrices $ D_i^k $, $ i\in \{1,\cdots,N\} $. 
However, since the recursion $ \hat{D}_{i}^k:= D_i^k  - M_U (\hat{D}_{i+1}^k)^{-1} M_L$ with $ \hat{D}_{N+1}^k :=0 $ needs be performed from $ i=N $ to $ 1 $ in the block Gaussian elimination method and $ \hat{D}_{i}^k $ is not necessarily sparse (dense matrix inversion  needs  to be performed), the lower-level sparsity is not preserved, which makes the structure-exploiting methods computationally expensive for a PDE system with a fine-grained spatial discretization, i.e., with a large $ n_x $.

\section{Double-layer Jacobi method}\label{sec_doublelayer}
This section introduces the proposed double-layer Jacobi method, which makes use of both the sparsities in the upper-level KKT matrix in \eqref{KKT_sys} and the lower-level Jacobian matrices $ D_i^k $, $ i\in \{1,\cdots,N\} $. 
The upper-layer Jacobi method is a general NMPC optimization method that is not limited to a particular class of dynamical systems. 
Convergence for the upper-layer Jacobi method is analyzed, and some variants of the Jacobi method are given. 
The lower-layer Jacobi method is a linear solver dedicated to the NMPC control of PDE systems by exploiting their particular structures after spatial discretization. 
We first give some general results on the convergence of iterative methods in the following subsection. 
\subsection{Preliminaries}
We first review the Jacobi method (see, e.g., \cite{saad2003iterative})  for solving linear equations as follows. 
\begin{lemma}\label{lemma_jaco}
	Let
	\begin{equation}\label{eq_linsys}
	Av = b,
	\end{equation}
	where $ A\in\mathbb{R}^{n\times n} $ and $ b\in\mathbb{R}^{n}$. 
	Let $ A $ be decomposed into $ A = D+L+U $, where $ D $, $ L $, and $ U $ are the diagonal, strict lower triangular, and strict upper triangular elements (blocks) of $ A $, respectively. 
	Assume that $ D $ is invertible. 
	The Jacobi method for solving \eqref{eq_linsys} is given by 
	\begin{equation*}
	v^{k+1} = D^{-1}(b-(L+U)v^k).
	\end{equation*}
	The Jacobi method converges if 
	\begin{equation}\label{eq_cvg_linsys}
	\rho(D^{-1}(L+U))<1.
	\end{equation}
	Likewise, for the element-wise decomposition, the Jacobi method converges if the matrix $ A $ is strictly diagonally dominant. 
\end{lemma}
We show in the following some results on the convergence of general iterations. 
\begin{definition}
	(Point of attraction \cite{ortega2000iterative}). Consider the iteration 	
	\begin{equation}\label{eq_poc}
	v^{k+1} = H(v^k),
	\end{equation}
	where $ v\in\mathbb{R}^{n}$ and $ H: P \rightarrow \mathbb{R}^{n} $ for a subset $ P \subset \mathbb{R}^{n} $. 
	Let $ v^*  $ be an interior point of $ P $ and a fixed point of the iteration \eqref{eq_poc}, i.e., $ v^{*} = H(v^*) $. 
	Then, $ v^* $ is said to be a point of attraction of the iteration \eqref{eq_poc} if there is an open neighborhood $ O\subset P $ of $ v^* $ such that the iterates defined by \eqref{eq_poc} all lie in $ O $ and converge to $ v^* $ for any $ v^0\in O $. 
\end{definition}
\begin{lemma}\label{lemma_ost}	 
	\cite{ortega2000iterative}
	Consider the iteration {\eqref{eq_poc}.  }
	Then, $ v^* $ is a point of attraction of the iteration \eqref{eq_poc} if the following condition holds:
	\begin{equation*}
	\rho(\nabla_v H(v^*))<1.
	\end{equation*}
\end{lemma}
\begin{lemma}\label{lemma_rate}
	(Convergence factor and rate \cite{ortega2000iterative}). If  $ v^*  $ is a point of attraction of the iteration \eqref{eq_poc}, the following holds:
	\begin{equation*}
	\rho(\nabla_v H(v^*)) = \lim_{k\rightarrow \infty}\sup\|v^k-v^*\|^{1/k},
	\end{equation*}
	and $ \rho(\nabla_v H(v^*)) $ is called the convergence factor. 
	The convergence rate is defined by $ -\ln \rho(\nabla_v H(v^*)) $.
\end{lemma}

\subsection{Upper-layer Jacobi method} 
Let $ D^k $, $ L $, and $ U $ be the diagonal, strict lower triangular, and strict upper triangular blocks of the KKT matrix in \eqref{KKT_sys} as follows.
\begin{equation*}
\begin{split}
&D^k:=\blkdiag(D_1^k,\cdots,D_N^k) \\
&L:= \text{lower-block-diag}(M_L,\cdots,M_L) \\
&U:= \text{upper-block-diag}(M_U,\cdots,M_U)
\end{split}
\end{equation*}
The block-diagonal matrix  $ D^k $ can be guaranteed to be invertible if $ h $ is sufficiently small and $ \gamma >0 $. 
The upper-layer Jacobi method for solving the KKT conditions \eqref{eq_KKT} is given by 
\begin{equation}\label{eq_1st_jacobi}
S^{k+1} = S^k - \alpha^{\max} (D^k)^{-1}\mathcal{K}^{k},
\end{equation}
where $ \alpha^{\max}\in (0,1] $ is a scalar obtained from the fraction-to-the-boundary rule \eqref{eq_alphamax} and $ S^0 $ is chosen such that the primal feasibility condition $ G(u_i,x_i)> 0 $ is satisfied for all $ i\in\{1,\cdots,N\} $.

The upper-layer Jacobi method exploits the banded structure (temporal sparsity) of the KKT matrix by ignoring its off-diagonal blocks. 
That is, the couplings introduced by the state and costate equations are ignored during iteration so that \eqref{eq_1st_jacobi} can be performed block-wisely ($ D^k $ is a block-diagonal matrix). 
Although the Jacobi method in \eqref{eq_1st_jacobi} can be regarded as Newton's method ignoring off-diagonal blocks and Newton's method is known to be locally quadratically convergent under mild assumptions, the convergent property might not be preserved for the Jacobi method. 
The convergence of the upper-layer Jacobi method is analyzed in the following subsection. 

\subsubsection{Convergence}
We first give a general convergence condition for the upper-layer Jacobi method.

\begin{thm}\label{theorem_cvg_rho} 
	$ S^* $ is a point of attraction of the iteration \eqref{eq_1st_jacobi} if the following condition holds:
	\begin{equation}\label{eq_jaco_cv_cdt}
	\rho((D^*)^{-1}(L+U))<1.
	\end{equation}
\end{thm}
\textbf{Proof.}
Since $ G(u_i^*,x_i^*)>0 $ is satisfied for all $ i\in \{1,\cdots,N\} $, there exists an  open neighborhood $ O $ of $ S^* $ such that for any $ S^0\in O $, the fraction-to-the-boundary rule \eqref{eq_alphamax} will never be triggered when the iteration converges, i.e., $ \alpha^{\max} = 1 $ in the neighborhood of  $ S^* $. 
The result then follows by applying Lemma \ref{lemma_ost} with $ \alpha^{\max} = 1 $, $ {\mathcal{K}}^{*} = 0 $, and $ \nabla_S{{\mathcal{K}}}^{*} = D^* +L+U$.
$\quad \Box $

Theorem \ref{theorem_cvg_rho} gives a general sufficient condition for convergence. However, the condition \eqref{eq_jaco_cv_cdt} can only be verified afterward and does not provide significant insights related to the NMPC problem. 
We show in the following lemma and theorem that the upper-layer Jacobi method can be guaranteed to converge by tuning the NMPC parameters, such as the prediction horizon $ T $ and the regularization parameter $ \gamma $. 

\begin{lemma}\label{lemma_rho0}
	Let $ D_i^*$ be decomposed  into 
	\begin{equation}\label{eq_decom}
	\begin{split}
	D_i^* 
	&= 
	\left[
	\begin{array}{ccc}
	h\nabla_{x}f^*_i- I & 0& 0  \\
	h\nabla_{ux}^2\mathcal{H}^*_i&  h\nabla_{uu}^2\mathcal{H}^*_i + \gamma I & 0\\
	h\nabla_{xx}^2\mathcal{H}^*_i & h\nabla_{xu}^2\mathcal{H}^*_i  & h(\nabla_{x}f^*_i)^T- I
	\end{array}
	\right] \\
	&\quad +
	h\left[
	\begin{array}{ccc}
	0 &  \nabla_{u}f^*_i & 0 \\
	0& 0 & (\nabla_{u}f^*_i)^T \\
	0&0 & 0
	\end{array}
	\right]\\
	&=:
	\bar{D}_i^* + h\tilde{D}_i^*.
	\end{split}
	\end{equation}
	Let $ \bar{D}^* $ and $ \tilde{D}^*  $ be defined as follows.
	\begin{equation*}
	\begin{split}
	&\bar{D}^*:=\blkdiag(\bar{D}_1^*,\cdots,\bar{D}_N^*)\\
	&\tilde{D}^*:=\blkdiag(\tilde{D}_1^*,\cdots,\tilde{D}_N^*)
	\end{split}
	\end{equation*}
	Then, for any $ h> 0 $  and $ \gamma\geq 0 $ such that $ \bar{D}^* $ is invertible, e.g., when $ h $ is sufficiently small  and $ \gamma  $ is nonzero, the following holds:
	\begin{equation*}
	\rho((\bar{D}^*)^{-1}(L+U))= 0.
	\end{equation*}
\end{lemma}
\textbf{Proof.}  See Appendix \ref{apx_rho0}.

\begin{thm}\label{theorem_cvg} 
	The convergence of the upper-layer Jacobi method \eqref{eq_1st_jacobi} is described as follows.
	\begin{enumerate}[(i)]
		\item Let $ \gamma  $ be chosen to be nonzero. There exists $ T_0>0 $ such that $ \rho((D^*)^{-1}(L+U))<1 $ holds for any $ T<T_0 $. \label{theorem_ii}
		\item There exists $ \epsilon > 0 $ such that $ \rho((D^*)^{-1}(L+U))<1 $ holds for any $ \tilde{D}_i^* $ satisfying $ \| \tilde{D}_i^* \| < \epsilon \|\bar{D}_i^*\|$, $ i\in\{1,\cdots,N\} $.  \label{theorem_iii}
	\end{enumerate}
\end{thm}
	\textbf{Proof.} 
	Proof of (\ref{theorem_ii}). 
	From the definitions in Lemma \ref{lemma_rho0}, we know that  
	\begin{equation}\label{eq_pertur}
	\rho((D^*)^{-1}(L+U)) = \rho((\bar{D}^*+h\tilde{D}^*)^{-1}(L+U)).
	\end{equation}
	We first choose $ T_0>0 $ to be sufficiently small. 
	Since $ \gamma  $ is a nonzero constant, $ T<T_0 $, and $ h=T/N $ is sufficiently small, the right-hand side of \eqref{eq_pertur} can be seen as a small perturbation of $ \rho((\bar{D}^*)^{-1}(L+U)) $ in terms of $ h $.
	That is, together with Lemma \ref{lemma_rho0}, we obtain that
	\begin{equation*}
	\lim_{h\rightarrow 0}\rho((D^*)^{-1}(L+U)) = 0.
	\end{equation*}
	From the continuities of matrix inverse and spectral radius, there always exists $T_0 > 0$ such that $ \rho((D^*)^{-1}(L+U))<1 $ holds for any $T < T_0$.

	Proof of (\ref{theorem_iii}). Let $ \epsilon>0 $ be chosen to be sufficiently small. 
	Since $ \epsilon $ is a small number and $ \| \tilde{D}_i^* \| < \epsilon \|\bar{D}_i^*\|$, $ i\in\{1,\cdots,N\} $, $ \rho((D^*)^{-1}(L+U)) = \rho((\bar{D}^*+h\tilde{D}^*)^{-1}(L+U)) $ can also be seen as a small perturbation of $ \rho((\bar{D}^*)^{-1}(L+U)) = 0 $ in terms of $ \tilde{D}^* $. 
	Similarly to the proof of (\ref{theorem_ii}), the result then follows. $\quad \Box $

Theorem \ref{theorem_cvg} can be  interpreted as follows. 
Theorem \ref{theorem_cvg} (\ref{theorem_ii}) indicates that a nonzero regularization parameter $ \gamma $ and a short prediction horizon $ T $ guarantee the convergence of the upper-layer Jacobi method. 
However, note that neither a nonzero $ \gamma $ nor a small $ T>0 $ is a necessary condition for satisfying $ \rho((D^*)^{-1}(L+U))<1 $. 
Since $ \bar{D}_i^* $ consists of the sensitivity $ \nabla_{u}f^*_i $, Theorem \ref{theorem_cvg} (\ref{theorem_iii}) can be interpreted as that  $ \rho((D^*)^{-1}(L+U))<1 $ holds if the dynamical system \eqref{eq_ode} is not sensitive to the control input $ u $. 

\begin{remark}\label{rmk_gamma}
	(Role of regularization).  
	As can be seen from the proof of Theorem \ref{theorem_cvg} (\ref{theorem_ii}), the nonzero regularization parameter $ \gamma $ makes $ (\bar{D}^*)^{-1} $ less ill-conditioned (less sensitive to $ h $) and therefore guarantees the convergence of the method when selecting a short prediction horizon $ T $. 
	However, a large $ \gamma $ makes the optimal solution of the RR-NMPC problem \eqref{eq:nlp_relaxed} far away from the original NMPC problem \eqref{eq:nlp} unless the regularization reference $ \tilde{u}_i^* $ in the regularization term \eqref{eq_reg} is chosen to be a good estimation of the optimal control input of \eqref{eq:nlp}. 
	Since the RR-NMPC problem has to be solved successively at every time step, $ \tilde{u}_i^* $ can be fixed to be the optimal control input of the last time step or updated at each iteration, i.e.,  $ \tilde{u}_i^* = u_i^k $.
	If $ \tilde{u}_i^* = u_i^k  $, the optimal solution to the RR-NMPC problem is equivalent to the optimal solution to the NMPC problem with only relaxation. 
\end{remark}

\subsubsection{Jacobi-type variants}\label{sec_jacovariants}
The reason the iteration \eqref{eq_1st_jacobi} is called the Jacobi method is that its convergence condition \eqref{eq_jaco_cv_cdt} mimics the condition \eqref{eq_cvg_linsys} for solving linear equations. 
This can be seen as applying the Jacobi method to solve the KKT conditions. 
Likewise, Jacobi-type methods (see, e.g., \cite{saad2003iterative}) for solving linear equations, such as the Gauss-Seidel and successive over-relaxation (SOR) methods, can be applied to solve the KKT conditions as well. 
We show several Jacobi-type methods and their conditions for convergence as follows.
\begin{itemize}
	\item Forward Gauss-Seidel method (FGS): 
	\begin{equation}\label{eq_fgs}
	S^{k+1} = S^k - \alpha^{\max}(D^k+L)^{-1}{\mathcal{K}}^{k}
	\end{equation}
	Condition for convergence:
	\begin{equation}\label{eq_cvg_fgs}
	\rho((D^*+L)^{-1}U)<1.
	\end{equation}
	\item Backward Gauss-Seidel method (BGS): 
	\begin{equation*}
	S^{k+1} = S^k - \alpha^{\max}(D^k+U)^{-1}{\mathcal{K}}^{k}
	\end{equation*}
	Condition for convergence:
	\begin{equation}\label{eq_cvg_bgs}
	\rho((D^*+U)^{-1}L)<1.
	\end{equation}
	\item SOR method:
	\begin{equation*}
	S^{k+1} = S^k - \alpha^{\max}\omega(D^k+\omega L)^{-1}{\mathcal{K}}^{k},\ \omega >0.
	\end{equation*}
	Condition for convergence:
	\begin{equation*}
	\rho((D^*+\omega L)^{-1}(\omega U+(\omega -1)D^*))<1.
	\end{equation*}
	\item Symmetric Gauss-Seidel (SGS) method: 
	\begin{equation}\label{eq_sgs}
	S^{k+1} = S^k - \alpha^{\max}(D^k+L)^{-1}({\mathcal{K}}^{k}-U(D^k+U)^{-1}{\mathcal{K}}^{k}).
	\end{equation}
	Condition for convergence:
	\begin{equation}\label{eq_sgs_cvg}
	\rho((D^*+ L)^{-1}U(D^*+ U)^{-1}L)<1.
	\end{equation}
\end{itemize}

Note that FGS and BGS have the same amount of computation as the Jacobi method. 
The difference is the rate of convergence as shown in the following theorem. 
\begin{thm}\label{cor_cvg}
	If the convergence condition \eqref{eq_jaco_cv_cdt} holds for the Jacobi method, then the convergence conditions \eqref{eq_cvg_fgs} and \eqref{eq_cvg_bgs} also hold for FGS and BGS, respectively. 
	That is, $ S^* $ is a point of attraction of the FGS and BGS iterations. 
	Moreover, both the FGS and BGS methods converge twice as fast as the Jacobi method. 
\end{thm}
\textbf{Proof.}
Since the KKT matrix in \eqref{KKT_sys} is a block-tridiagonal matrix, the KKT matrix is consistently ordered \cite{hageman2012applied}. 
It is known from \cite{saad2003iterative} that for a consistently ordered matrix, the spectral radius of FGS is the square of that of the Jacobi method, i.e.,
\begin{equation*}
\rho((D^*+L)^{-1}U)  = \rho((D^*)^{-1}(L+U))^2.
\end{equation*}
The conclusion above can be shown similarly for BGS that 
\begin{equation*}
\rho((D^*+U)^{-1}L) = \rho((D^*)^{-1}(L+U))^2.
\end{equation*}
Recall the definition of the convergence rate in Lemma \ref{lemma_rate}. 
The result then follows. $\quad \Box $


As can be seen from  Theorems \ref{theorem_cvg} and \ref{cor_cvg}, a short prediction horizon $ T $  and a nonzero $ \gamma  $ can also guarantee the convergence of FGS and BGS. 
The discussion on the role of the regularization procedure in Remark \ref{rmk_gamma} applies to FGS and BGS as well. 
As for the SGS iteration \eqref{eq_sgs}, it can be seen as a BGS iteration followed by a FGS iteration, i.e., a backward sweep followed by a forward sweep. 
The iteration in the previous work \cite{zavala2016new} is similar to the FGS iteration \eqref{eq_fgs}. 
However, the inequality constraints are kept and the regularization procedure is not introduced, so the convergence is difficult to guarantee.

\subsection{Lower-layer Jacobi method}
The upper-layer Jacobi method and its variants essentially consist of solving linear equations with the coefficient matrices $ D_i^k $ of the structure in \eqref{eq_decom} for $ i\in\{1,\cdots,N\} $. 
Since $ D_i^k $ is sparse and its structure is fixed, efficient exact or iterative solution methods usually exist.  
For example, the Jacobi method for solving linear equations can be applied directly when $ h $ is sufficiently small. 
In this subsection, we introduce another iterative method by exploiting the particular structure of $ D_i^k $.

The linear systems are reordered to have the following coefficient matrix:
\begin{equation}\label{eq_sublinsys_cpx}
\left[
\begin{array}{lll}
0 & h\nabla_{x}f^k_i -  I & h\nabla_{u}f^k_i \\
h(\nabla_{x}f^k_i)^T - I& h\nabla_{xx}^2\mathcal{H}^k_i & h\nabla_{xu}^2\mathcal{H}^k_i\\
h(\nabla_{u}f^k_i)^T& h\nabla_{ux}^2\mathcal{H}^k_i  &  h\nabla_{uu}^2\mathcal{H}^k_i + \gamma I
\end{array}
\right].
\end{equation}
For the sake of brevity, the linear system with the coefficient matrix \eqref{eq_sublinsys_cpx} is expressed by using the following shorthand:
\begin{equation}\label{eq_sublinsys_simple}
\left[
\begin{array}{lll}
0 & F_x& F_u \\
F_x^T& A_{xx} & A_{xu}\\
F_u^T& A_{ux}  &  A_{uu}
\end{array}
\right]
\left[
\begin{array}{c}
v_{1}\\ 
v_2\\ 
v_3
\end{array} 
\right] = 
\left[
\begin{array}{c}
b_1\\ 
b_{2}\\ 
b_3
\end{array} 
\right].
\end{equation}
Equation \eqref{eq_sublinsys_simple} can be solved by first eliminating $ (v_{1},v_2) $ and then solving for $ v_3 $, i.e., by performing the following two steps. 
\begin{subequations}\label{eq_step12}
\begin{equation}\label{eq_step1}
\begin{split}
\left(
A_{uu}
- 
\left[
\begin{array}{cc}
F_u^T&A_{ux} 
\end{array} 
\right]
\left[
\begin{array}{cc}
0 & F_x \\ 
F_x^T & A_{xx}
\end{array} 
\right]^{-1}
\left[
\begin{array}{c} 
F_u \\  
A_{xu} 
\end{array} 
\right]
\right) 
v_3\\
=
b_3 - 
\left[
\begin{array}{cc}
F_u^T&A_{ux} 
\end{array} 
\right]
\left[
\begin{array}{cc}
0 & F_x \\ 
F_x^T & A_{xx}
\end{array} 
\right]^{-1}
\left[
\begin{array}{c} 
b_1 \\  
b_{2} 
\end{array} 
\right]
\end{split}
\end{equation}
\begin{equation}\label{eq_step2}
\begin{split}
\left[
\begin{array}{cc}
0 & F_x \\ 
F_x^T & A_{xx}
\end{array} 
\right]
\left[
\begin{array}{c} 
v_{1}  \\  
v_2
\end{array} 
\right] =
\left[
\begin{array}{c} 
b_1 \\  
b_{2} 
\end{array} 
\right]
-
\left[
\begin{array}{c} 
F_u \\  
A_{xu} 
\end{array} 
\right]
v_3
\end{split}
\end{equation}
\end{subequations}
Solving \eqref{eq_step12} consists of solving several linear equations of the following form:
\begin{equation}\label{eq_coff}
\left[
\begin{array}{cc}
0 & F_x \\ 
F_x^T & A_{xx}
\end{array} 
\right]\left[
\begin{array}{c} 
v_4 \\  
v_5 
\end{array} 
\right]=\left[
\begin{array}{c} 
b_4 \\  
b_5 
\end{array} 
\right].
\end{equation}
The linear equation \eqref{eq_coff} can be solved by first solving $ F_xv_5= b_4 $ and then solving $ F_x^Tv_4 = b_5-A_{xx}v_5  $. 
That is, linear equations with the coefficient matrices $ F_x $ and $ F_x^T $ are solved essentially. 
We show next that these linear equations can be solved efficiently by using the Jacobi method.

Recall that for the discretized PDE system \eqref{eq_ode}, we have  
\begin{equation*}
F_x = h\nabla_{x}f^k_i - I = 
\left[
\begin{array}{cc}
-I & hI \\ 
h\nabla_W g^k_i & h\nabla_{\dot{W}} g^k_i-I \\ 
\end{array} \right].
\end{equation*} 
Therefore, a linear equation with the coefficient matrix $ F_x $, i.e., the following equation,
\begin{equation*}
\left[
\begin{array}{cc}
-I & hI \\ 
h\nabla_W g^k_i & h\nabla_{\dot{W}} g^k_i-I \\ 
\end{array} \right]\left[
\begin{array}{c} 
v_6 \\  
v_7 
\end{array} 
\right]=\left[
\begin{array}{c} 
b_6 \\  
b_7 
\end{array} 
\right]
\end{equation*}
can be solved by first eliminating $ v_6 $ with $ v_6=hv_7-b_6 $ and then solving a linear system with the following coefficient matrix:
\begin{equation}\label{eq_jaco2_sys}
h\nabla_{\dot{W}} g^k_i-I + h^2\nabla_W g^k_i.
\end{equation}
Since that $ \nabla_{\dot{W}} g^k_i $ is diagonal, the off-diagonal entries of \eqref{eq_jaco2_sys} have orders of $ \mathcal{O}(h^2) $.
It is easy to show that \eqref{eq_jaco2_sys} is diagonally dominant if $ h $ is sufficiently small. 
Moreover, the off-diagonal entries of $ \nabla_W g^k_i $ are sufficiently small for  many of the PDE equations, such as the heat transfer equation and the Navier-Stokes equation with a large Reynolds number.
%
That is, according to Lemma \ref{lemma_jaco}, convergence of the Jacobi method can be achieved for solving linear equations with the coefficient matrices \eqref{eq_jaco2_sys}. 
The same conclusion can be made for the $ F_x^T $ system. 

The lower-layer Jacobi method is concluded as follows. Since $ n_x \gg n_u $ for PDE-constrained NMPC problems, the major computational cost for solving \eqref{eq_sublinsys_simple} comes from solving linear equations with the coefficient matrices \eqref{eq_jaco2_sys}, which can be solved efficiently by using the Jacobi method if, e.g., the NMPC problem is finely discretized in time, i.e.,  with a sufficiently small $ h $.

\begin{remark}\label{rmk_second_general}
	The key point of the lower-layer method is to find a method that solves the equation \eqref{eq_sublinsys_simple} efficiently, e.g., the introduced lower-layer Jacobi method for the NMPC control of the PDE system \eqref{eq_pde}. 
	For PDE systems that are not in the form of \eqref{eq_pde} or even general systems, the proposed upper layer's iteration can be performed efficiently if a structure-exploiting linear solver can be designed. 
\end{remark} 

\section{Numerical experiment}\label{sec_exp}
In this section, we demonstrate the performance of the double-layer Jacobi method in terms of the computation time, number of iterations, and convergence factor by using a heat transfer closed-loop control example. 
The experiment was implemented in C and performed on a 3.9-GHz (turbo boost frequency) Intel Core i5-8265U laptop computer. 
To reduce the effect of the computing environment, the computation time at each time step was measured by taking the minimum one of ten runs of the closed-loop simulation. 
\subsection{System description}
We consider a nonlinear heat transfer process in a thin copper plate \cite{matlab2dheat}. 
Because the plate is relatively thin compared with the planar dimensions, temperature can be assumed constant in the thickness direction. 
The system is described by the following two-dimensional PDE:
\begin{equation*}
\rho C_p t_z \frac{\partial w(p,t)}{\partial t} - k t_z \triangle w(p,t) +2 Q_c + 2 Q_r = 0,
\end{equation*}
where $ w $ is the plate temperature, $ p\in \Omega:=\{(x,y)|x,y\in[0,1]\} $ ($ x $ here stands for the horizontal axis) and $ Q_c $ and $ Q_r $ are, respectively, the convection and radiation heat transfers defined as follows.
\begin{equation*}
\begin{split}
&Q_c := h_c(w(p,t)-T_a) \\
&Q_r := \epsilon \delta (w(p,t)^4-T_a^4)
\end{split}
\end{equation*}
The boundary conditions are the zero Neumann boundary conditions. 
The parameters of the heat transfer process are given in Table \ref{tab_para}. 
\begin{table}[!ht]
	\caption{Parameters in the heat transfer process}
	\scalebox{0.87}{\begin{tabular}{l|l|l}
			\hline 
			$ \rho $ & 8960  & Density of copper [kgm$^{-3}$] \\ 
			\hline 
			$ C_p $ & 386  & Specific heat of copper [Jkg$ ^{-1}$K$ ^{-1}$] \\ 
			\hline 
			$ t_z $ & 0.01  & Plate thickness [m] \\ 
			\hline 
			$ k $ & 400  & Thermal conductivity of copper [Wm$ ^{-1}$K$ ^{-1} $] \\ 
			\hline 
			$ h_c $ & 1  & Convection coefficient [Wm$ ^{-2}$K$ ^{-1} $] \\ 
			\hline 
			$ T_a $ & 300  & Ambient temperature [K] \\ 
			\hline 
			$ \epsilon $ & 0.5  & Emissivity of the plate surface \\ 
			\hline 
			$ \delta $ & $ 5.67\cdot 10^{-8} $  & Stefan-Boltzmann constant [Wm$ ^{-2}$K$ ^{-4}$] \\ 
			\hline 
	\end{tabular}}
	\label{tab_para}
\end{table}
There are 16 actuators distributed uniformly under the plate to heat or cool the plate above the ambient temperature in the range of $ [T_a, T_a+400] $ K.  
The temperatures  of the plate at the positions of the actuators can be controlled directly.
We assume that the actuators negligibly impact the convection and radiation heat transfer processes. 
We are interested in controlling the temperature distribution across the plate under constraints, which is a typical control problem that arises in semiconductor manufacturing. For example, a temperature gradient needs to be maintained within a wafer to ensure catalytic activation \cite{bleris2006towards}. 
\subsection{NMPC description}
\begin{figure}
	\begin{center}
		\includegraphics[width=0.3\textwidth]{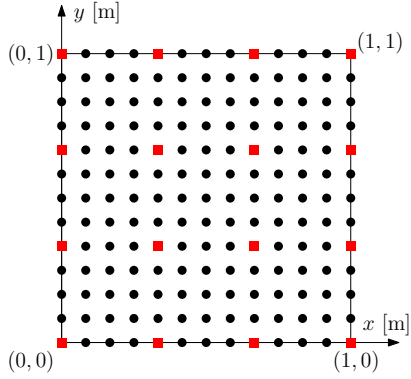}
	\end{center}
	\caption{Spatial discretization grid points on plate (red squared points: actuators' positions)}\label{fig_2ddis}
\end{figure}

The plate was uniformly discretized into $ 13\times 13 $ spatial grid points as shown in Fig. \ref{fig_2ddis}.
Since the temperatures at the positions of the actuators can be controlled directly, the temperatures of the red squared points are regard as control inputs. 
We obtain a system with 16 inputs and 153 states. 
The inputs are constrained by 
\begin{equation*}
G(u,x) = 
\left[
\begin{array}{c}
u - T_ae\\ 
-u +(T_a+400) e 
\end{array} 
\right]\geq 0,
\end{equation*}
where $ e =[1,\cdots,1]^T $. 
We chose the cost function to be quadratic as
\begin{equation*}
l_i(u,x):=\frac{1}{2}(\|x-x_{{ref}}\|^2_{Q} + \|u-u_{{ref}}\|^2_{R}),\  i\in\{1,\cdots,N\},
\end{equation*}
where $ x_{{ref}} $ and $ u_{ref} $ encoded the temperature distribution reference and the weighting matrices were $ Q=I $ and $ R= 0.1\times I $. 

Note that the lower-layer Jacobi method converged fast and needed only two iterations due to the small thermal diffusivity $ k\rho^{-1}C_p^{-1} $.
Moreover, we noticed that the coefficient matrix in \eqref{eq_step1} was dominated by the diagonal matrix $ A_{uu} $. 
The linear equation \eqref{eq_step1} in the lower-layer Jacobi method was solved iteratively by performing two of the following iterations:
\begin{equation*}
\begin{split}
&A_{uu}v_3^{k+1}
=
b_3 +\\ 
&\left[
\begin{array}{cc}
F_u^T&A_{ux} 
\end{array} 
\right]
\left[
\begin{array}{cc}
0 & F_x \\ 
F_x^T & A_{xx}
\end{array} 
\right]^{-1}
\left(
\left[
\begin{array}{c} 
F_u \\  
A_{xu} 
\end{array} 
\right]
v_3^k
-
\left[
\begin{array}{c} 
b_1 \\  
b_2 
\end{array} 
\right]
\right),
\end{split}
\end{equation*}
which can be solved efficiently due to the diagonal property of $ A_{uu} $. 
The parameters of the NMPC controller are given in Table \ref{tab_NMPC_para}. 
\begin{table}[!ht]
	\caption{NMPC parameters}
	\scalebox{0.87}{\begin{tabular}{l|l}
			\hline 
			Name & Value  \\ 
			\hline 
			Prediction horizon $ T $ & 100 [s]  \\ 
			\hline 
			\# of temporal discretization points $ N $ & 20   \\ 
			\hline 
			Barrier parameter $ \tau $ & 100 \\ 
			\hline
			Regularization reference $ \tilde{u}_i^* $ & $ u_i^k $ \\ 
			\hline 
			Regularization parameter $ \gamma $ & 0.5 \\ 
			\hline 
			Stopping criterion & $ \|\mathcal{K}^{k}\|_{\infty}<1 $\\ 
			\hline 
			Upper-layer method & SGS \eqref{eq_sgs}    \\ 
			\hline 
			Lower-layer method & Two Jacobi iterations \\ 
			\hline 
	\end{tabular}}
	\label{tab_NMPC_para}
\end{table}
Since all of the matrices during iteration were sparse, the expressions of the matrix-vector multiplications were pre-computed offline, which made the proposed method matrix-free. 
The computational complexity of the proposed method for the heat transfer example is $ \mathcal{O}(N(n_x+n_u)) $.

\subsection{Closed-loop simulation}
\begin{figure}[!ht]
	\centering
	\subfloat[$ 0\leq t\leq 500 $ s]{\includegraphics[width = 0.24\textwidth]{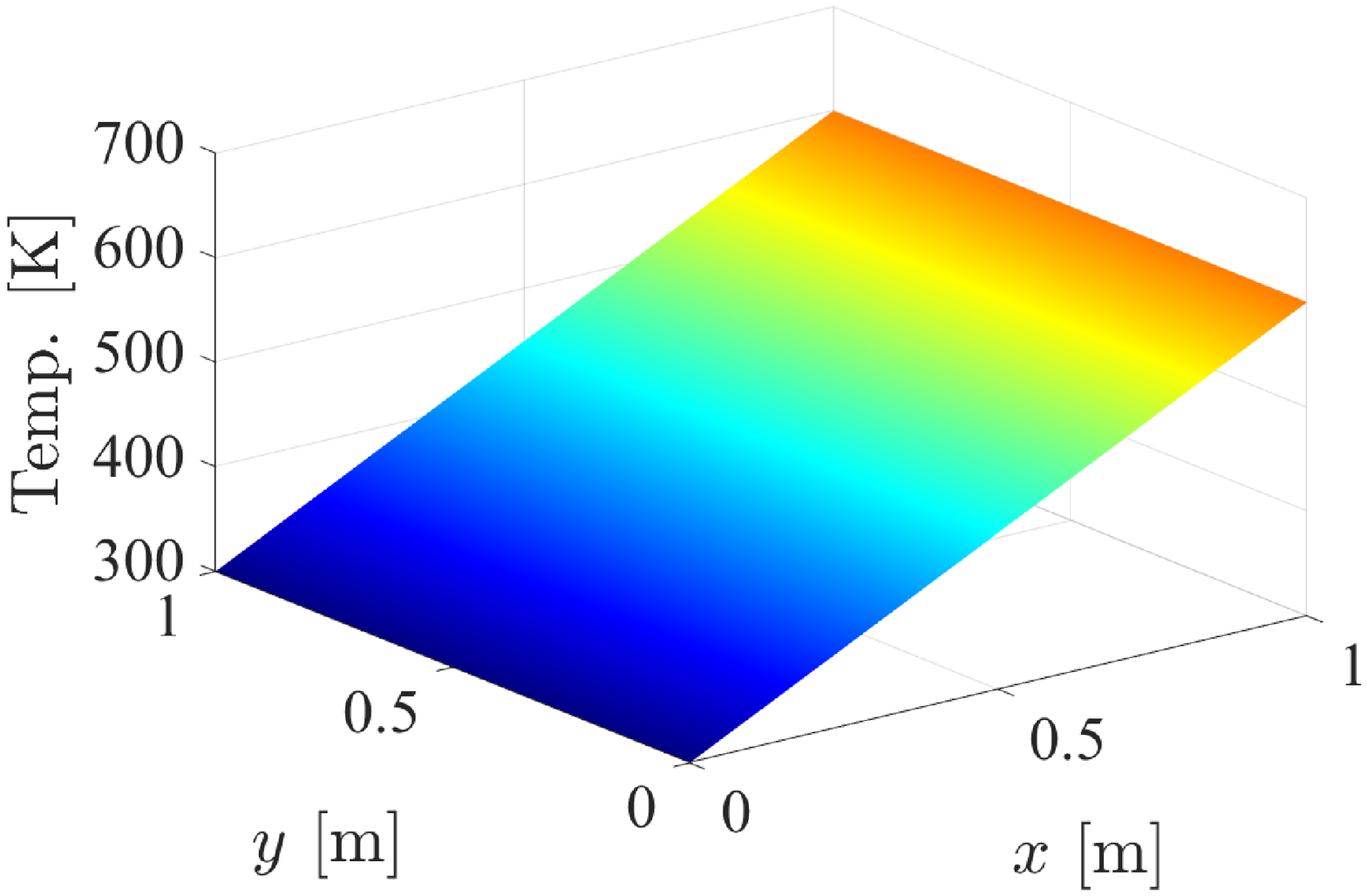}}
	\subfloat[$ 500< t\leq 1000 $ s]{\includegraphics[width = 0.24\textwidth]{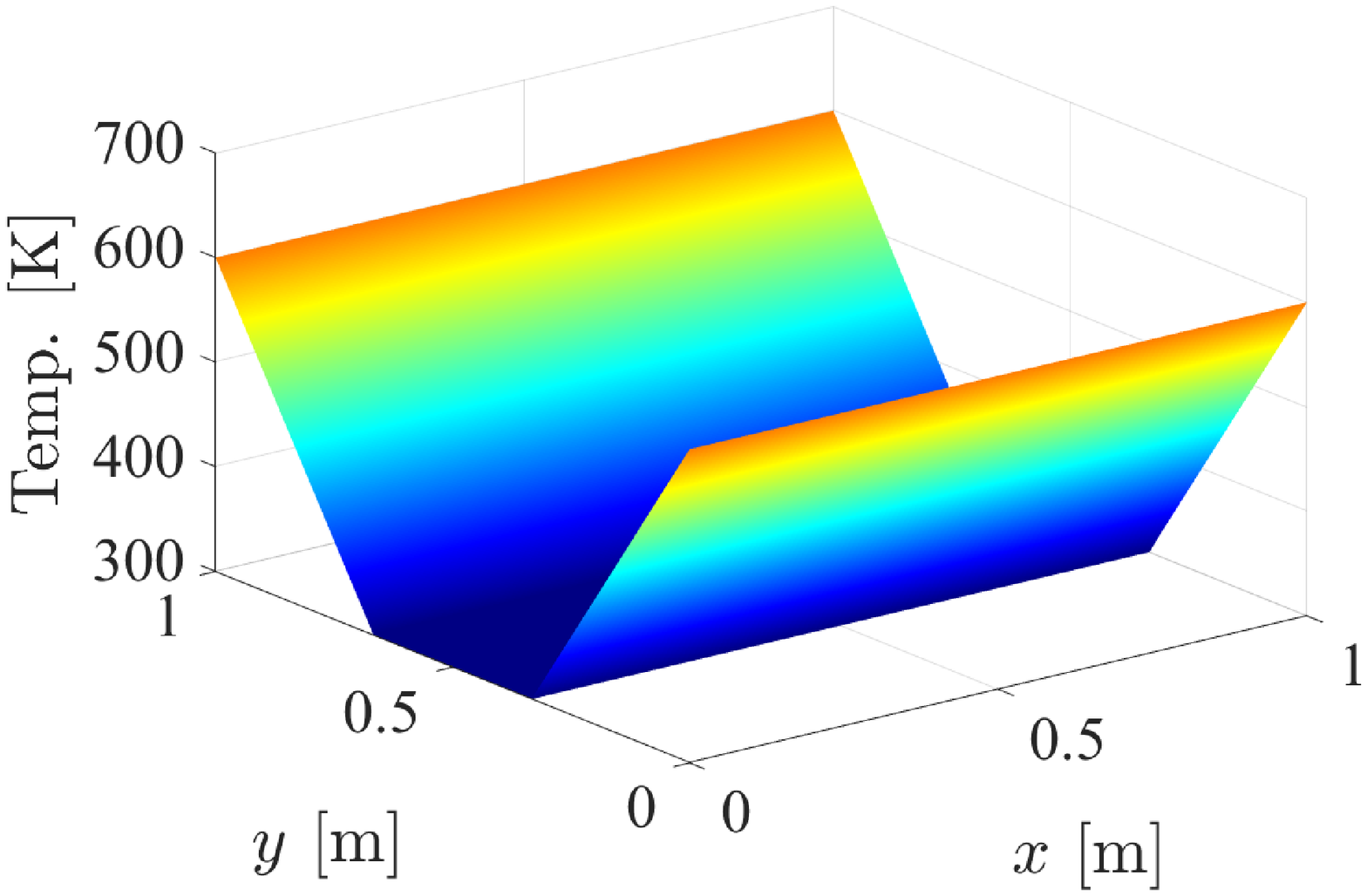}}
	\caption{Temperature distribution references at different time periods}
	\label{fig_2dref}
\end{figure}

The system was started from an initial state of $ \bar{x}_0 =[T_a,\cdots,T_a]^T $. 
The simulation was performed for $ 1000 $ s with a sampling period of $ 5 $ s. 
The temperature distribution reference, as shown in Fig. \ref{fig_2dref}, was set to a slope shape for the first 500 seconds and a V-like shape for the last 500 seconds. 
The second reference was fed to the controller by changing continuously from the first reference within 50 seconds. 

\begin{figure}[!ht]
	\centering
	\subfloat[$ t=50 $ s]{\includegraphics[width = 0.24\textwidth]{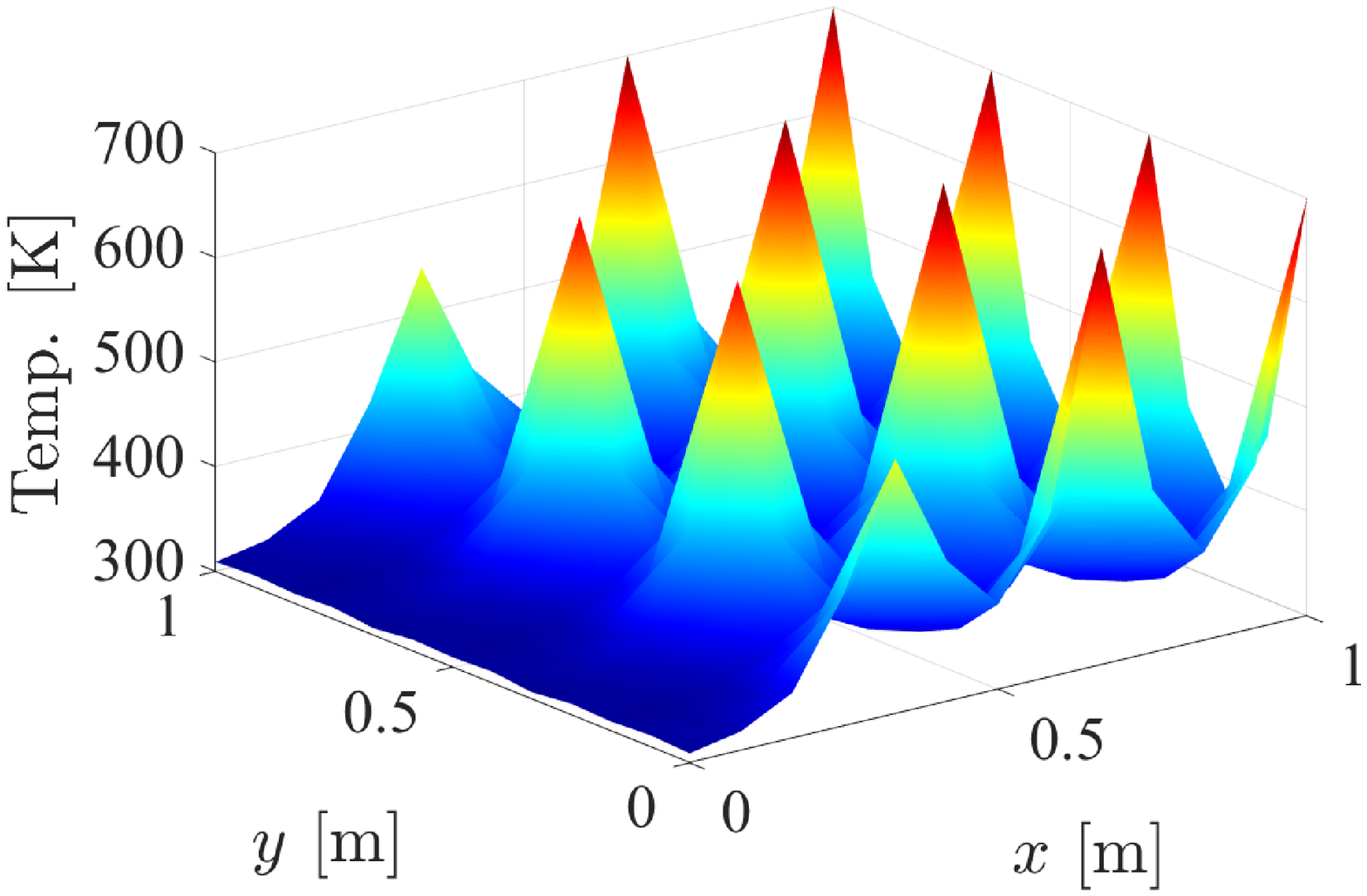}} 
	\subfloat[$ t=500 $ s]{\includegraphics[width = 0.24\textwidth]{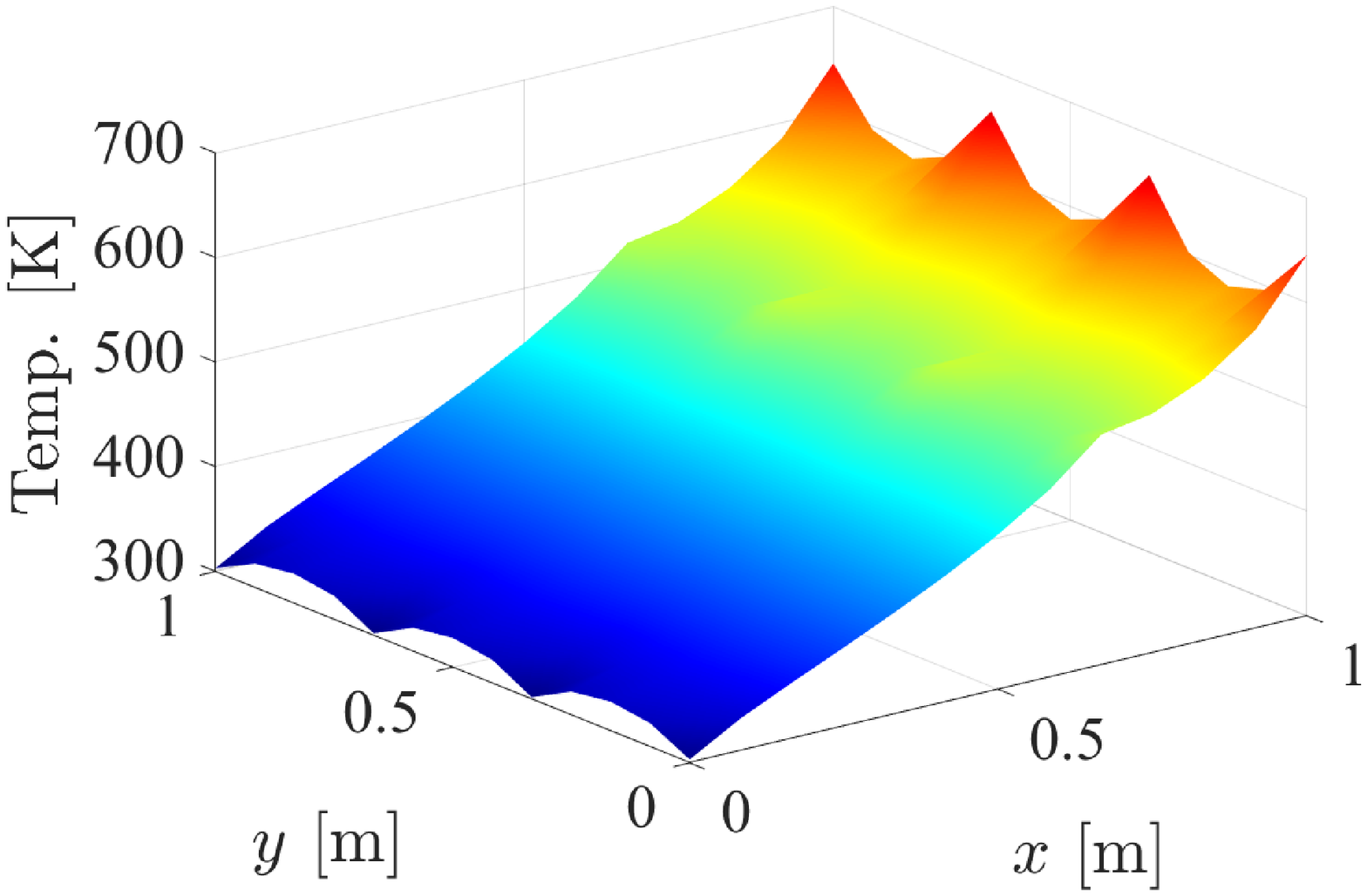}}\\
	\subfloat[$ t=550 $ s]{\includegraphics[width = 0.24\textwidth]{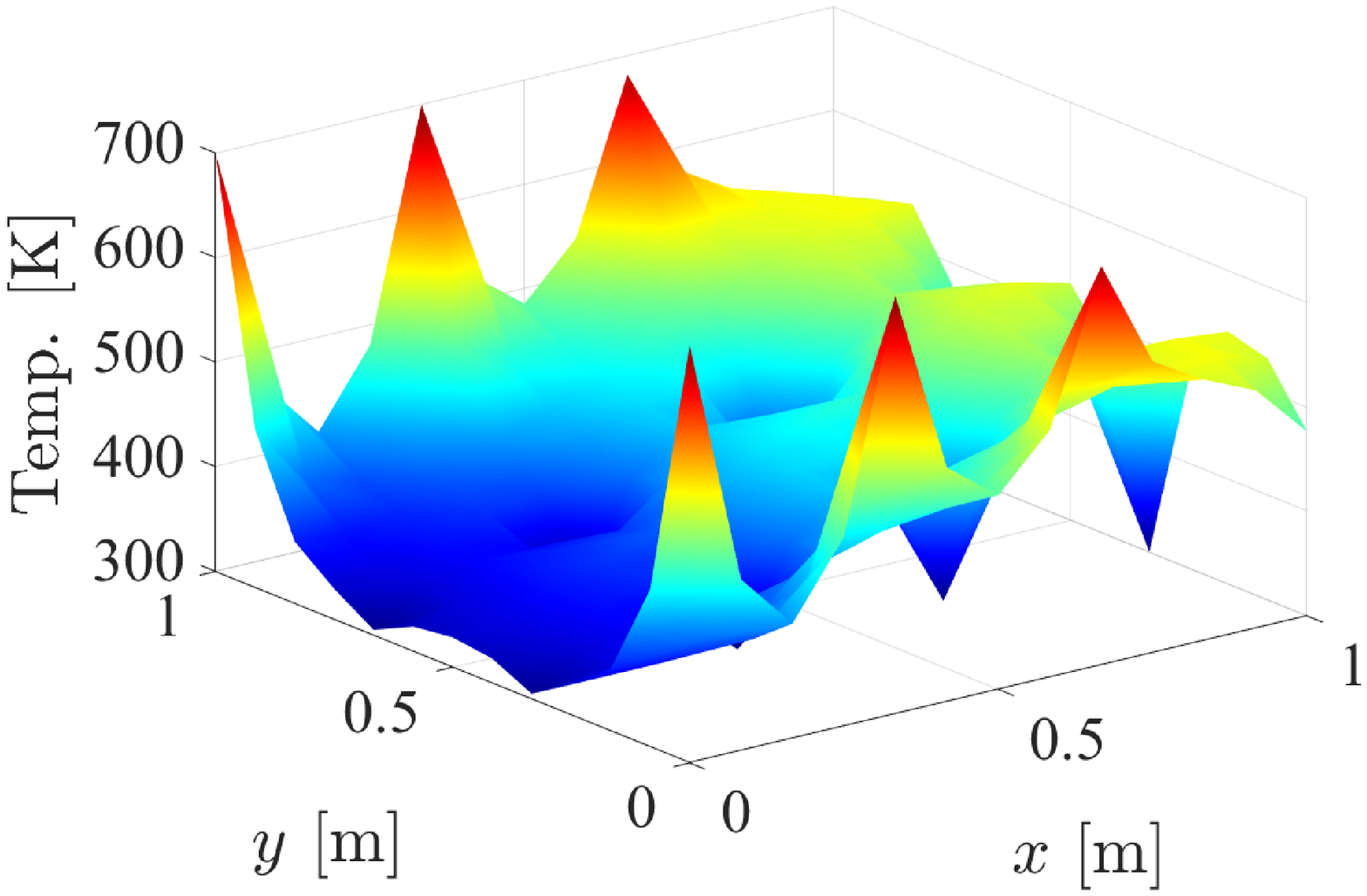}}
	\subfloat[$ t=1000 $ s]{\includegraphics[width = 0.24\textwidth]{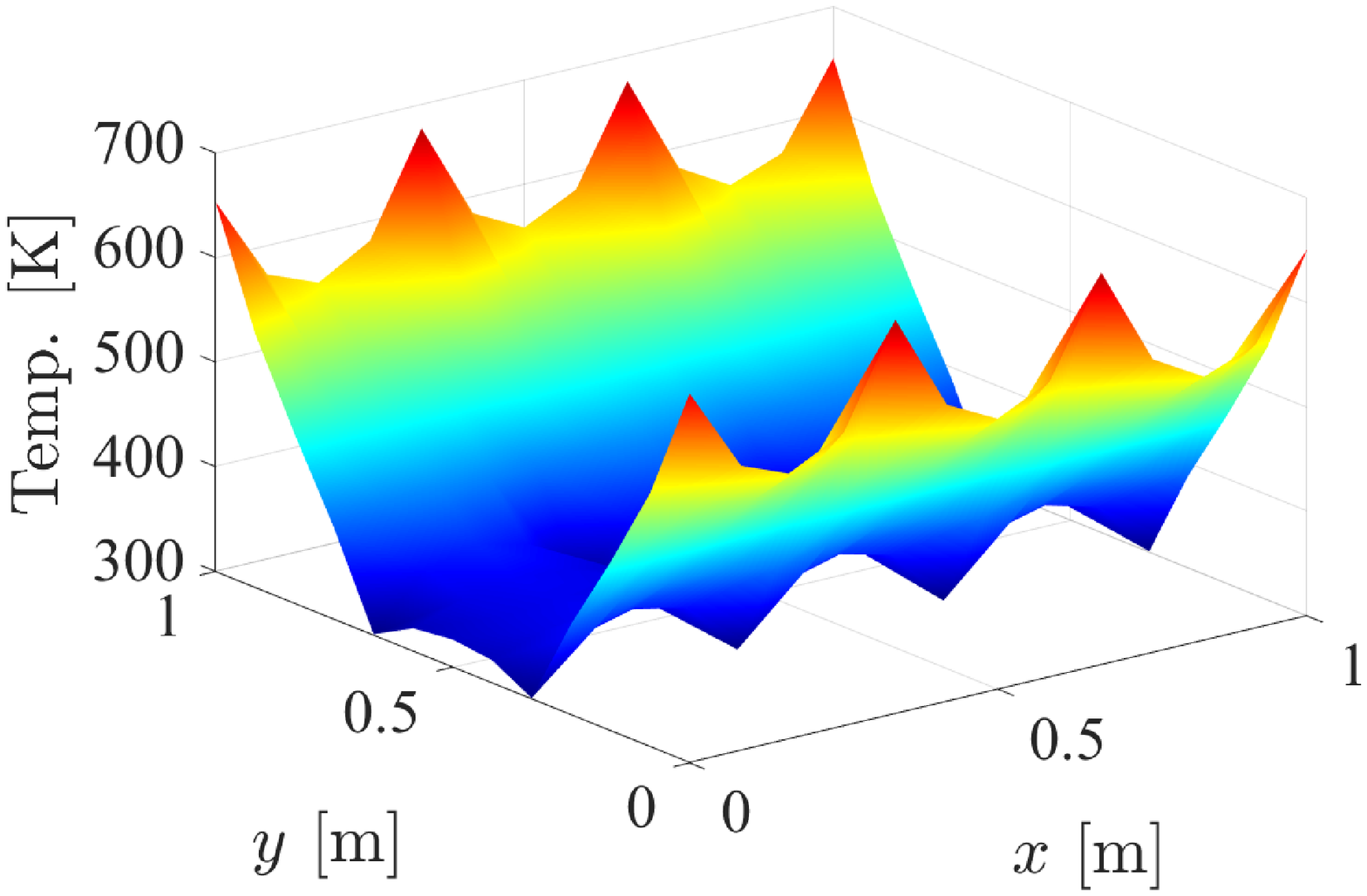}} 
	\caption{Temperature distributions at different time $ t $}
	\label{fig_2dres}
\end{figure}

\begin{figure}[!htb]
	\begin{center}
		\includegraphics[width=0.5\textwidth]{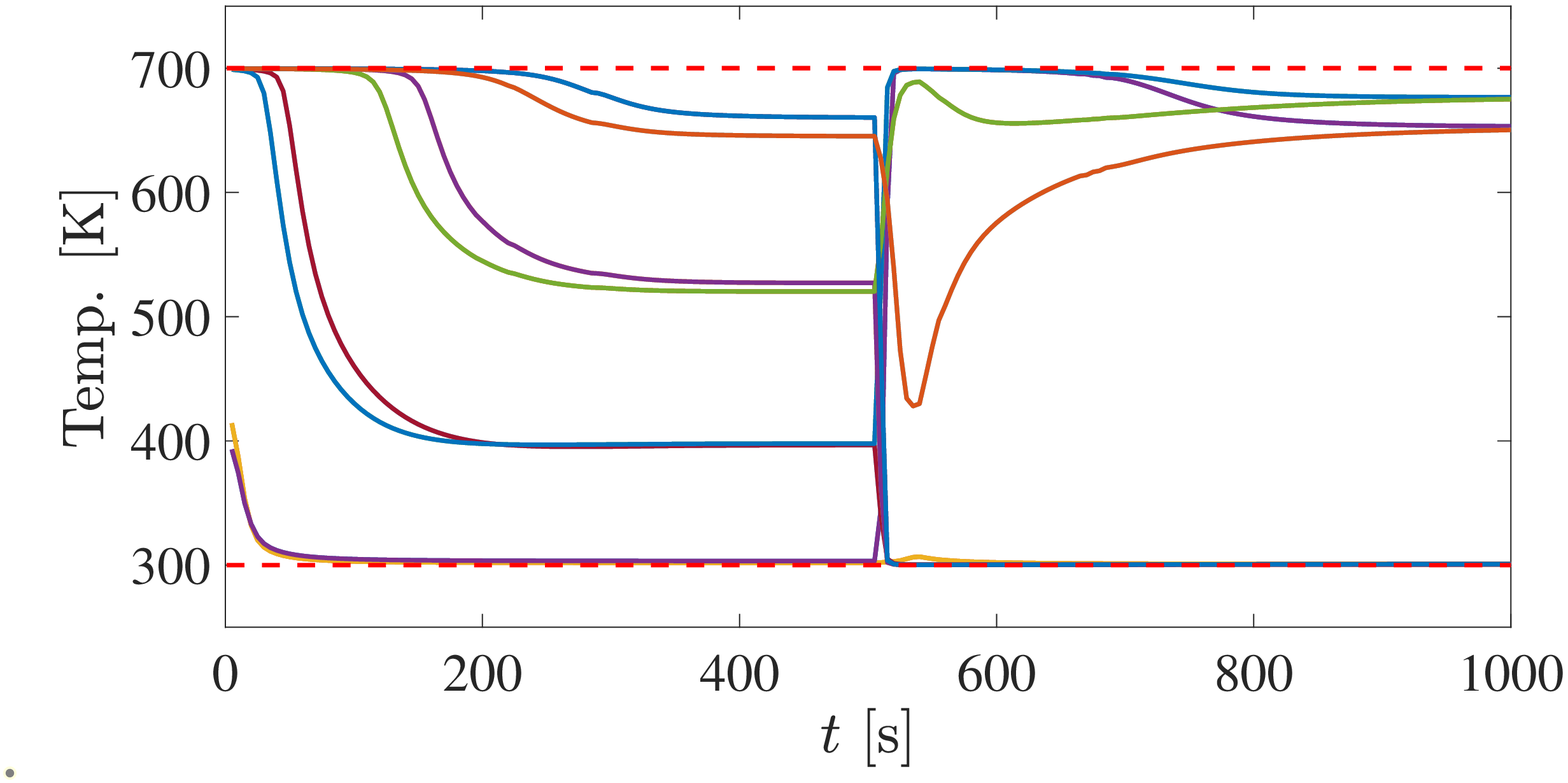}
	\end{center}
	\caption{Time histories of inputs (some inputs coincide with each other)}\label{fig_2dresu}
\end{figure}

For tracking the first reference of the closed-loop simulation, two sampled plots at $ t=50 $ s and $ t=500 $ s are shown in Fig. \ref{fig_2dres} (a) and (b). 
For the second reference, two sampled plots at $ t=550 $ s and $ t=1000 $ s are shown in Fig. \ref{fig_2dres} (c) and (d). 
As shown by these plots, the references were tracked well by using the NMPC controller. 
The time histories of the control inputs are shown in Fig. \ref{fig_2dresu}.  
Although the barrier parameter $ \tau  $ was fixed to $ 100 $, a high accuracy was still achieved such that the inputs approached the boundaries very closely. 

\begin{figure}[!htb]
	\begin{center}
		\includegraphics[width=0.5\textwidth]{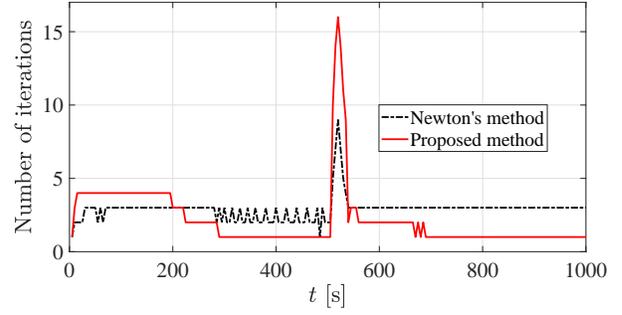}
	\end{center}
	\caption{Time histories of numbers of iterations}\label{fig_2dnoi}
\end{figure}

\begin{figure}[!htb]
	\begin{center}
		\includegraphics[width=0.5\textwidth]{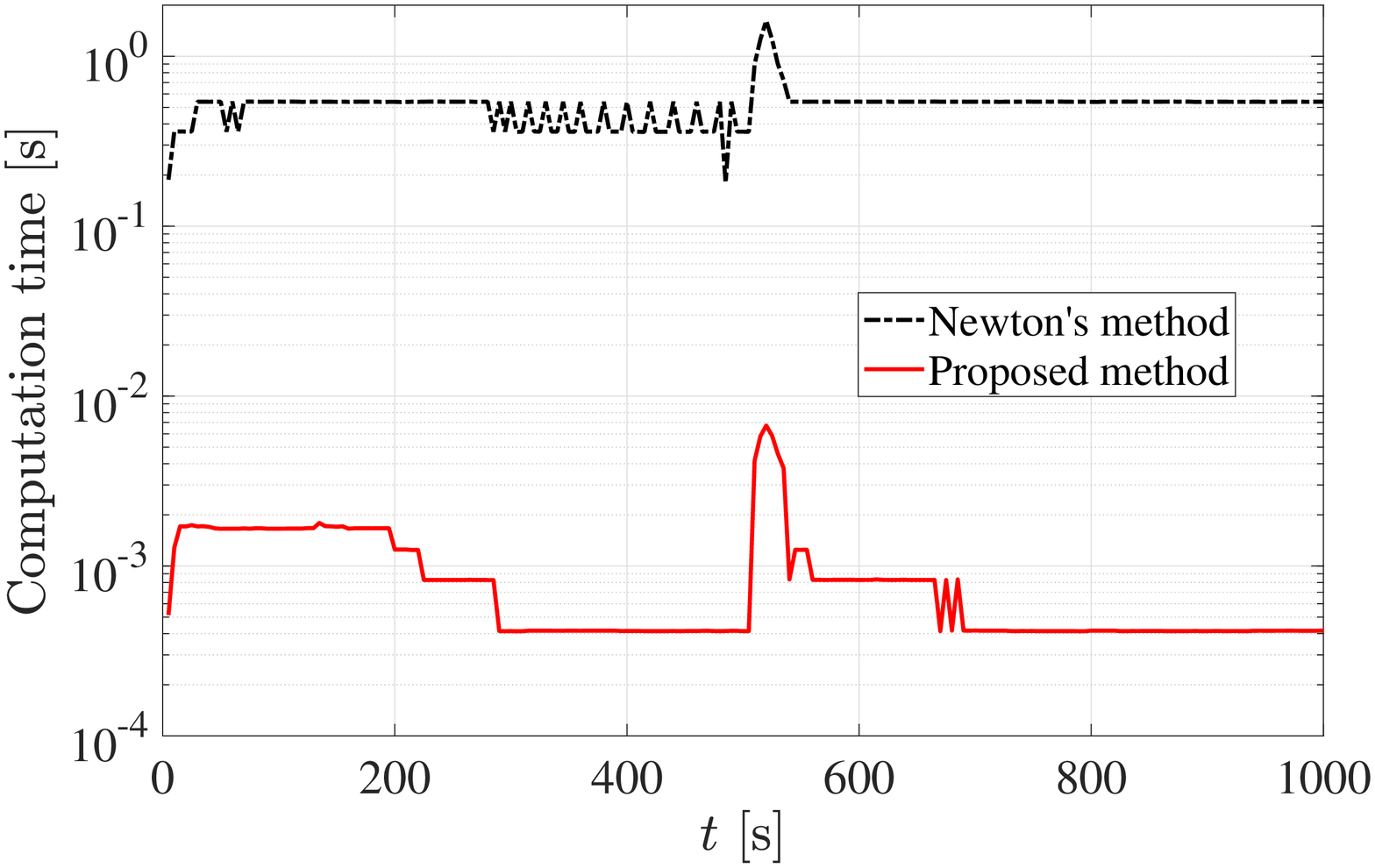}
	\end{center}
	\caption{Time histories of computation times per time step}\label{fig_2dct}
\end{figure}

To demonstrate the performance of the proposed method, we compared the conventional Newton's method introduced in Section \ref{sec_reg_newton}. 
The search direction \eqref{KKT_sys} in Newton's method  was calculated by using the block Gaussian elimination method, which was implemented by using NMPC real-time optimization software ParNMPC \cite{deng2018parallel}. 
Note that since both methods were based on the interior-point method, their computation times per iteration were consistent throughout the closed-loop simulation. 
The mean computation time per iteration for Newton's method was 0.180 s, which was about $ 433 $ times of that of the proposed method (0.416 ms). 
Considering that their numbers of iterations shown in  Fig. \ref{fig_2dnoi} were in the same range, the proposed method was much faster than Newton's method in terms of the computation time per time step shown in Fig. \ref{fig_2dct}. 
Furthermore, since the proposed method for the heat transfer example was matrix-free, its compiled executable file size (440 KB) was only about one tenth of that of Newton's method, which enables embedded applications for the proposed method. 

\begin{figure}[!htb]
	\begin{center}
		\includegraphics[width=0.5\textwidth]{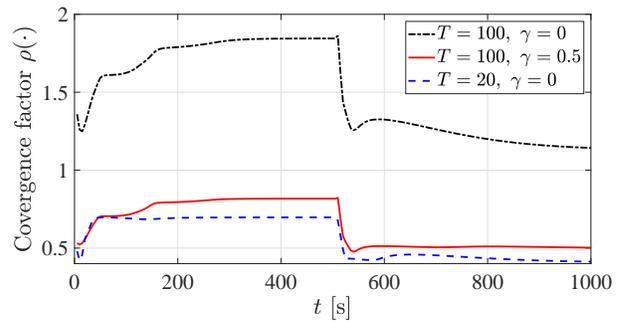}
	\end{center}
	\caption{Time histories of convergence factors of SGS under different settings}\label{fig_2dcvg}
\end{figure}

Lastly, we discuss the effects of the prediction horizon and regularization. 
In the numerical experiment, the proposed method could not converge without regularization ($ \gamma=0 $). 
According to Theorem \ref{theorem_cvg} and Remark \ref{rmk_gamma}, the convergence can be guaranteed by shortening the prediction horizon and introducing a positive regularization parameter $ \gamma $. 
We compared the convergence factor $ \rho((D^*+ L)^{-1}U(D^*+ U)^{-1}L) $ for the upper layer's SGS iteration along the closed-loop simulation under different prediction horizons ($ T=20$ and $100 $) and regularization parameters ($ \gamma =0$ and $0.5 $) in Fig. \ref{fig_2dcvg}. 
It can be seen that the convergence condition $ \rho((D^*+ L)^{-1}U(D^*+ U)^{-1}L)<1 $ was satisfied with either regularization or a short prediction horizon.

\section{Conclusion}\label{sec_conclusion}
This paper presents a double-layer Jacobi method for the NMPC control of PDE systems. 
The NMPC problem is formulated on the basis of the spatially and temporally discretized PDE system and then relaxed and regularized. 
The proposed method performs simple Jacobi-type iterations to solve the KKT conditions and the underlying linear systems to make full use of the sparsities exist in both the spatial and temporal directions. 
Furthermore, the convergence of the proposed method can be guaranteed by adjusting the prediction horizon and regularization parameter.
The results of the numerical experiment show that the proposed method can significantly reduce the computation time and program size.

Future research directions include extending the proposed method to the NMPC control of other large-scale systems and applying the finite element method to discretize the PDE system.

\appendix
\section{Proof of Lemma \ref{lemma_rho0}}\label{apx_rho0}
If $ N=1 $ {($ L=U=0 $)}, the result can be easily obtained. We discuss the case of $ N\geq 2 $. 
The proof is done by showing that the eigenvalues of $ (\bar{D}^*)^{-1}(L+U) $ are all zero. 
{In fact,} the expression  $ \det (\sigma I - (\bar{D}^*)^{-1}(L+U)) = \sigma^{N(2n_x+n_u)}$ is obtained by using Schur complement recursively {as shown below.} 

Define a set of matrices 
\begin{equation*}
\mathbb{A}
:=
\left\{
\left[
\begin{array}{ccc}
0 & 0 & 0 \\ 
0 & 0 & 0 \\ 
P & 0 & 0
\end{array} 
\right]\in \mathbb{C}^{(2n_x+n_u)\times (2n_x+n_u)},\ P\in \mathbb{C}^{n_x\times n_x}
\right\}.
\end{equation*}
Define the following shorthand:
\begin{equation*}
\bar{D}_i^L: = -(\bar{D}_i^*)^{-1}M_L\ \text{and}\ \bar{D}_i^U: = -(\bar{D}_i^*)^{-1}M_U\ 
\end{equation*}
so that 
	\begin{equation*}
	(\bar{D}^*)^{-1}(L+U) 
	=
	-\left[
	\begin{array}{ccccc}
	0 & \bar{D}_{1}^U   &        &        &          \\
	\bar{D}_{2}^L  & 0 & \bar{D}_{2}^U    &        &          \\
	& \bar{D}_{3}^L   & \ddots & \ddots &          \\
	&        & \ddots & 0 & \bar{D}_{N-1}^U     \\
	&        &        & \bar{D}_{N}^L   & 0
	\end{array}
	\right].
	\end{equation*}
	For any $ A\in \mathbb{A} $, $ \sigma\in \mathbb{C} $, and $ i\in\{2,\cdots,N\} $, it can be examined that 
	\begin{equation}\label{eq_inA}
	\bar{D}_{i-1}^U(\sigma I-A)^{-1}\bar{D}_{i}^L \in \mathbb{A}.
	\end{equation}
Let $ K\in \{2,\cdots,N\} $ and $ A_K\in \mathbb{A}  $. 
We define a $K $-size ($ K $ blocks of rows and columns) block-tridiagonal matrix $ W_K $ by
\begin{equation}\label{eq_mblock}
\begin{split}
\left[
\begin{array}{cccc|c}
\sigma I & \bar{D}_{1}^U   &        &        &          \\
\bar{D}_{2}^L  & \sigma I & \bar{D}_{2}^U    &        &          \\
& \bar{D}_{3}^L   & \ddots & \ddots &          \\
&        & \ddots & \sigma I & \bar{D}_{K-1}^U     \\
\hline
&        &        & \bar{D}_{K}^L   & \sigma I-A_{K}
\end{array}
\right]
=:
\left[
\begin{array}{cc}
M_A & M_B \\ 
M_C & M_D
\end{array} 
\right].
\end{split}
\end{equation}
The determinant of $ W_K $ is given by 
\begin{equation}\label{eq_det_schur}
\det W_K = \det (\sigma I - A_K) \det (\text{schur}(W_K,\sigma I - A_K)),
\end{equation}
where $ \text{schur}(W_K,\sigma I - A_K) $ denotes the Schur complement of the block $ \sigma I - A_K $ of $ W_K $, i.e.,
\begin{equation*}
\text{schur}(W_K,\sigma I - A_K) = M_A - M_BM_D^{-1}M_C.
\end{equation*}
Let us then calculate the right hand side of \eqref{eq_det_schur}. 
It can be shown that 
\begin{equation}\label{eq_eig_sigma}
\det (\sigma I - A_K)  = \sigma^{2n_x+n_u}.
\end{equation}
Since $ A_K\in \mathbb{A}  $, we can know from \eqref{eq_inA} that the only nonzero block (lower right  corner) of $ M_BM_D^{-1}M_C $ belongs to $\mathbb{A}  $, i.e,
\begin{equation}\label{eq_schur}
\bar{D}_{K-1}^U(\sigma I-A_K)^{-1}\bar{D}_{K}^L \in \mathbb{A}.
\end{equation}
By choosing $ A_{K-1}  $ to be the left hand side of \eqref{eq_schur}, the Schur complement $ \text{schur}(W_K,\sigma I - A_K) $ can be seen as a $( K-1) $-size block-tridiagonal matrix in the form of \eqref{eq_mblock}, i.e.,
\begin{equation}\label{eq_schur_recur}
\text{schur}(W_K,\sigma I - A_K) =: W_{K-1}.
\end{equation}
By substituting \eqref{eq_eig_sigma} and \eqref{eq_schur_recur} into \eqref{eq_det_schur}, we obtain the following recursion:
\begin{equation*}
\det W_K = \sigma^{2n_x+n_u} \det W_{K-1}. 
\end{equation*}
Following the procedures above and together with  {$  W_1 =  \sigma I_{{2n_x+n_u}}  $}, we obtain
\begin{equation*}
\det W_K = \sigma^{K(2n_x+n_u)},
\end{equation*}
which holds for any $ K\in \{2,\cdots,N\} $ and $ A_K\in\mathbb{A} $. 
Then, by choosing  {$ K=N $ and $ A_N=0\in\mathbb{A} $}, we have 
\begin{equation}\label{eq_wm}
\det W_N = \det(\sigma I - (\bar{D}^*)^{-1}(L+U)) = \sigma^{N(2n_x+n_u)}.
\end{equation}
From \eqref{eq_wm}, we can know that $ (\bar{D}^*)^{-1}(L+U) $ has only zero eigenvalues. The conclusion $ \rho((\bar{D}^*)^{-1}(L+U)) =0 $ then follows.
$\quad \Box $

\bibliographystyle{IEEETrans}        
\bibliography{references}           

\end{document}